\documentclass[11pt]{article}
\usepackage{amscd}
\usepackage{amsfonts}
\usepackage{amsmath}
\usepackage{amssymb}
\usepackage{amsthm}
\usepackage{bbm}
\usepackage{CJK}
\usepackage{fancyhdr}
\usepackage{graphicx}
\usepackage{hyperref}
\usepackage{indentfirst}
\usepackage{latexsym}
\usepackage{mathrsfs}
\usepackage{xypic}

\newtheorem{thm}{Theorem}[section]
\newtheorem{lemma}[thm]{Lemma}
\newtheorem{dfn}[thm]{Definition}
\newtheorem{prop}[thm]{Proposition}
\newtheorem{exm}[thm]{Example}

\newcommand{\prf}{\mdseries \textbf{\textsc{Proof.}}\quad}

\usepackage[top=1in,bottom=1in,left=1.25in,right=1.25in]{geometry}
\def\<{\langle}
\def\>{\rangle}

\def\c{\cdot}

\def\D{\Delta}

\def\g{\gamma}

\def\o{\otimes}

\def\v{\varepsilon}

\date{}
\begin{document}
\renewcommand{\baselinestretch}{1.2}
\renewcommand{\arraystretch}{1.0}
\title{\bf Enveloping actions and duality theorems for  partial twisted  smash products }
\author {{\bf Shuangjian Guo$^{1}$, Shengxiang Wang$^{2}$\footnote{Correspondence: Shengxiang Wang, E-mail: wangsx-math@163.com}}\\
{\small 1. School of Mathematics and Statistics, Guizhou University of Finance and Economics}\\
{\small Guiyang  550025, P. R. China}\\
{\small 2. School of Mathematics and Finance, Chuzhou University}\\
 {\small Chuzhou 239000, P. R. China}}
 \maketitle
\begin{center}
\begin{minipage}{12.cm}

{\bf Abstract}  In this paper, we first generalize the theorem about
the existence of an enveloping action to a partial twisted  smash product. Second, we construct a Morita context between the partial twisted smash product and the twisted  smash
product related to the enveloping action.  Furthermore, we
show some results relating partial actions and partial representations over the partial twisted smash products,  which generalize the results of
Alves and Batista (Comm. Algebra, 38(8): 2872-2902, 2010). Finally, we  present versions of the duality theorems
of  Blattner-Montgomery for partial twisted smash products.
\smallskip

\smallskip

 {\bf AMS(2010) Subject Classification}  16T05; 16T15
 \end{minipage}
 \end{center}
 \normalsize

\section{ Introduction}
\def\theequation{0. \arabic{equation}}
\setcounter{equation} {0}

Partial group actions were considered first by Exel  in the context of operator
algebras and they turned out to be a powerful tool in the study of $C^{*}$-algebras
generated by partial isometries on a Hilbert space in \cite{E94}. A treatment from a purely algebraic point of view was given
recently in \cite{D07}, \cite{DE05}, \cite{DE00}. In particular, the algebraic study of partial actions and partial representations was initiated in \cite{DE05} and \cite{DE00}, motivating investigations in diverse directions. Now, the results are formulated in a purely algebraic way
independent of the $C^{*}$-algebraic techniques which originated them.

The concepts of partial actions and partial coactions of Hopf algebras on algebras were introduced by Caenepeel and
Janssen in \cite{C08}. In which they put the Galois theory for partial
group actions on rings into a broader context, namely, the partial entwining structures. In particular, partial actions of a group
 $G$ determine partial actions of the group algebra $kG$ in a natural way.  Further developments in the theory of partial Hopf actions were done by  Lomp in \cite{L08}.

   Alves and  Batista  extended
 several results from the theory of partial group actions to the Hopf algebra setting, they constructed
 a Morita context relating the fixed point subalgebra for partial actions of finite dimensional Hopf
 algebras, and constructed the partial smash product in \cite{A09}. Later, they constructed a Morita context between the partial smash
 product and the smash product related to the enveloping action, defined partial representations of Hopf algebras and
showed some results relating partial actions and partial representations in \cite{A10}.

Furthermore, they proved a dual version of the
globalization theorem: Every partial coaction of a Hopf algebra admits
an enveloping coaction. they explored some consequences
of globalization theorems in order to present versions of the duality theorems
of Cohen-Montgomery and Blattner-Montgomery for partial Hopf actions in \cite{A11}. Recently, they introduced partial representations
of Hopf algebras and gave the paradigmatic examples of them, namely, the partial
representation defined from a partial action and the partial representation related
to the partial smash product in \cite{A13}.

 Based on the existing results, this paper is organized as follows:

In Section 3, we introduce the notion of partial actions of a Hopf algebra containing
a partial  left action and a partial right action, and define a partial bimodule of a Hopf algebra.
Then we introduce a new notion-partial twisted smash product $\underline{A\circledast H}$, generalizing the twisted smash product in \cite{W98}.
Furthermore, we prove the  existence of an enveloping action for such a
partial twisted smash product.

In Section 4,  we construct a Morita context between the partial twisted smash
 product $\underline{A\circledast H}$ and the twisted smash product $B\circledast H$,
where $H$ is a Hopf algebra which acts partially on the unital
algebra $A$, $B$ is an enveloping action for partial actions. This
result can also be found in \cite{A10} for the context of partial group actions.

In Section 5,  we show that, under some certain conditions on the algebra $A$, the partial twisted smash
product $\underline{A\circledast H}$ carries a partial representation of $H$.

In Section 6, we explore some consequences of globalization theorems in order to present versions of the duality theorems
of  Blattner-Montgomery for partial twisted smash products.

 \section{Preliminaries}
\def\theequation{\arabic{section}.\arabic{equation}}
\setcounter{equation} {0}

Throughout the paper, let $k$ be a fixed
 field and all algebraic systems are supposed to be over $k$.
 Let $M$ be a vector space
 over $k$ and let $id_M$  the usual identity map.
  For the comultiplication
 $\D $ in a coalgebra $C$ with a counit $\v _C$,
 we use the Sweedler-Heyneman's notation (see Sweedler \cite{S69}):
 $ \Delta(c)=c_{(1)}\o c_{(2)}$, for any $c\in C$.
 \vskip0.2cm

 We first recall some basic results and propositions that we will need later from Alves and  Batista \cite{A09},\cite{A10}.
 \medskip

  \noindent{\bf 2.1.  Partial left module algebra}~ Let $H$ be a Hopf algerba  and $A$ an  algebra. $A$ is said to be a partial
   left $H$-module  algebra if there exists a $k$-linear map $\rightharpoonup=\{\rightharpoonup: H\otimes A\rightarrow A\}$ satisfying the following conditions:
\begin{eqnarray*}
&&h\rightharpoonup (ab) =(h_{(1)})\rightharpoonup a)(h_{(2)}\rightharpoonup b),\\
&& 1_H\rightharpoonup a = a,\\
&&  h\rightharpoonup(g \rightharpoonup a) =(h_{(1)}\rightharpoonup1_A)(h_{(2)}g\rightharpoonup a),
\end{eqnarray*}
for all $h,g\in H$ and $a, b\in A$.

\section{Enveloping actions}
\def\theequation{\arabic{section}.\arabic{equation}}
\setcounter{equation} {0}

In the context of partial actions of Hopf algebras, it is proved that a partial action of a
Hopf algebra on a unital algebra $A$ admits an enveloping action $(B, \theta)$ if and only if each of
the ideals $\theta(A)\unrhd B$ is a unital algebra in \cite{A10}.
In this section, we mainly extend this famous result to partial twisted smash products.
\medskip

 Now, we give the definition of a partial right $H$-module algebra similar to \cite{A09} as follows:
\begin{dfn}
   Let $ H $ be a Hopf algebra and $A$  an  algebra, $A$ is said to be a partial
   right $H$-module algebra if there exists a  $k$-linear map $\leftharpoonup=\{\leftharpoonup: A\otimes H\rightarrow A\}$ satisfying the following conditions:
  \begin{eqnarray*}
&&(ab)\leftharpoonup h =( a\leftharpoonup h_{(1)})( b\leftharpoonup h_{(2)}),\\
&& a\leftharpoonup 1_H = a,\\
&& (a\leftharpoonup g)\leftharpoonup h =(1_A\leftharpoonup h_{(1)})(a\leftharpoonup gh_{(2)}),
\end{eqnarray*}
for all $h,g\in H$ and $a, b\in A$.
\end{dfn}
\begin{dfn}
Let $H$ be a Hopf algebra  and $A$  an algebra. $A$ is
called a partial $H$-bimodule algebra if the following conditions hold:

(i) $A$ is not only a partial left $H$-module algebra with the partial left module action
$¡±\rightharpoonup¡±$ but also a partial right $H$-module algebra with the partial right module action
$¡±\leftharpoonup¡±$.

(ii) These two partial module
structure maps satisfy the compatibility condition, i.e., $(h\rightharpoonup a) \leftharpoonup g =h\rightharpoonup (a \leftharpoonup g)$
for all $a\in A$ and $h,g \in H$.

\end{dfn}

Let $H$ be a Hopf algebra with an antipode $S$ and $A$ a partial $H$-bimodule algebra.
We first propose a multiplication on the vector space $A\o H$:
$$(a\circledast h)(b\circledast g)=a(h_{(1)}\rightharpoonup b \leftharpoonup S(h_{(3)}))\circledast h_{(2)}g,$$
for all $a,c\in A$ and  $g,h\in H$. It is obvious that the multiplication is associative. In order to make it to be an unital algebra, we project onto the
$$\underline{A\circledast H}=(A\otimes H)(1_A\otimes 1_{H}).$$
Then we can deduce directly the form and the properties of typical
elements of this algebra $$\underline{a \circledast h}= a(h_{(1)}\rightharpoonup 1_A \leftharpoonup S(h_{(3)}))\otimes h_{(2)},$$   and finally verify that the product among typical elements satisfy
\begin{eqnarray}
(\underline{a\circledast h})(\underline{b\circledast g})=\underline{a(h_{(1)}\rightharpoonup b \leftharpoonup S(h_{(3)}))\circledast h_{(2)}g},
 \end{eqnarray}
for all $h,g\in H$ and $a, b\in A$.
\medskip

From the above definitions, we have

\begin{prop}
With the  notations as above, $ \underline{A\circledast H}$ is an associative algebra
with a multiplication given by Eq.(3.1) and with the unit $\underline{1_A \circledast 1_{H}}$, and call it by a partial twisted smash product,
  where $1_A$ is the unit of $A$.
\end{prop}
\prf Similar to \cite{A09}.
$\hfill \Box$

\begin{dfn}
Let $H$ and $A$ be Hopf algebras.  A skew pair is a triple $(A,H,\sigma)$ endowed with a
$k$-linear maps $\sigma: A\otimes H\rightarrow k$ such that the following conditions are satisfied.

(1) $\sigma(ab,h)=\sigma(a,h_{(1)}))\sigma(b,h_{(2)}),$

(2) $\sigma(a_{(1)},h)\sigma(a_{(2)},g)= \sigma(1_A,g_{(1)})\sigma(a,g_{(2)}h)= \sigma(1_A,h_{(1)})\sigma(a,gh_{(2)}),$

(3) $\sigma(a,1)=\varepsilon(a),$

\noindent for all $h, g\in H$ and $a, b\in A$.
\end{dfn}

\begin{exm}
Let $H$ be a Hopf algebra with a bijective antipode $S$  and $A$  a Hopf algebra.
Suppose that $(A,H,\sigma)$ is a skew pair, then we can define two actions of $H$ and $A$: for any $h\in H, b\in A$,
\begin{eqnarray*}
 h \rightharpoonup b &=&  b_{(2)}\sigma(b_{(1)}, h), \\
 b\leftharpoonup h &=&  b_{(1)}\sigma(b_{(2)},(S^{-1})^2(h)).
 \end{eqnarray*}
It follows that
\begin{eqnarray*}
 \underline{a\circledast  h}&=&a((h_{(1)}\rightharpoonup1_A\leftharpoonup S(h_{(3)}))\otimes h_{(2)}\\
 &=&\sigma(1_A,h_{(1)}))a\otimes h_{(2)}\sigma(1_A,S^{-1}(h_{(3)})).
 \end{eqnarray*}
It is not hard to verify that $(A,\rightharpoonup,\leftharpoonup)$ is a partial $H$-bimodule  algebra and the multiplication of $\underline{A\circledast H}$  is
\begin{eqnarray*}
(\underline{a\circledast  h})(\underline{b\circledast  k})
=\sigma(b_{(1)},(h_{(1)}))\underline{ab_{(2)}\circledast  h_{(2)}k}\sigma(b_{(3)},S^{-1}(h_{(3)})),
\end{eqnarray*}
for all $h,k\in H$ and  $a, b\in A$.
Then $\underline{A\circledast_{\sigma}H}$ is a partial twisted smash product.
\end{exm}

\begin{exm}
 As a $k$-algebra, the four dimensional Hopf algebra $H_{4}$ is generated by
two symbols $c$ and $x$ which satisfy the relations $c^{2}=1$,
 $x^{2}=0$ and $xc+cx=0$. The coalgebra structure on $H_4$ is determined by
$$
\Delta(c)=c\otimes c,\ \
\Delta(x)=x\otimes 1
+c\otimes x, \ \v(c)=1, \ \v(x)=0.
$$

Consequently, $H_{4}$ has the basis l (identity), $c,\ \ x, \ \
 cx$, we now consider the dual $H_{4}^{\ast}$ of $H_4$.
 We have $H_{4}\cong H_{4}^{\ast}$ (as Hopf algebras) via
 $$
 1\mapsto 1^{\ast}+c^{\ast},\ \ c\mapsto 1^{\ast}-c^{\ast},\ \ x\mapsto x^{\ast}+(cx)^{\ast},\ \ cx\mapsto x^{\ast}-(cx)^{\ast}.
 $$
Here $\{1^{\ast}, c^{\ast}, x^{\ast}, (cx)^{\ast}\}$ denote the dual
 basis of $\{1,c,x,cx\}$. Let $T=1^{\ast}-c^{\ast},~
 P=x^{\ast}+(cx)^{\ast},~TP=x^{\ast}-(cx)^{\ast}$, we
 get another basis $\{1,T,P,TP \}$ of $ H_{4}^{\ast}$.
 Recall from \cite{C08} if $A$ is the subalgebra $k[x]$ of $H_4$, it is shown that $A$ is a right partial $H_4$-comodule algebra with the coaction $$\rho(1)=\frac{1}{2}(1\o 1+1\o c+1\o cx), ~\rho^r(x)=\frac{1}{2}(x\o 1+x\o c+x\o cx).$$
 By similar way we can show that $A$ is a left partial $H_4$-comodule algebra with the coaction
 $$\rho(1)=\frac{1}{2}(1\o 1+c\o 1+cx\o 1), ~~\rho^l(x)=\frac{1}{2}(1\o x+c\o x+cx\o x).$$
It is clear that  $A$ is a partial $H_4$-bicomodule algebra. Then $A$ is a  partial $H_4^{\ast}$-bimodule algebra via
\begin{eqnarray*}
f\rightharpoonup a=\sum <f,a_{[1]}>a_{[0]},~a\leftharpoonup g=<g,a_{[-1]}>a_{[0]},~a\in A,f,g\in H^{\ast}.
 \end{eqnarray*}
 Therefore we can obtain the partial twisted smash product $k[x] \circledast H_4^{\ast}$, we only consider the elements $P, T$ of $H_4^{\ast}$ as follows, then
 \begin{eqnarray*}
 (x\o T)(x\o P)&&=x(T_{(1)}\rightharpoonup x\leftharpoonup S^{\ast}(T_{(3)})) \o T_{(2)}P\\
 &&=\sum x(T\rightharpoonup x\leftharpoonup S^{\ast}(T)) \o TP\\
 &&=\sum x<T, \frac{1}{2}(1+c+cx)>x<T, \frac{1}{2}(1+c+cx)> \o TP=0.
 \end{eqnarray*}
\end{exm}
\begin{dfn}
Let $H$ be a Hopf algebra and $A, B$ be two partial $H$-bimodule  algebras.
A morphism of algebras $\theta:A\rightarrow B$ is said to be a  morphism of  partial $H$-bimodule algebras if
 $\theta(h\rightharpoonup a\leftharpoonup k)=h\rightharpoonup \theta(a)\leftharpoonup k$ for all $h,k\in H$ and
 $a\in A$.  If, in addition,  $\theta$ is an isomorphism, the partial actions are called equivalent.
 \end{dfn}
\begin{lemma}
Let $H$ be a Hopf algebra, $B$ a $H$-bimodule  algebra and  $A$  an ideal of $B$ with
unity $1_A.$  Then $H$ acts partially on $A$ by
$h\rightharpoonup a = 1_A(h\triangleright a),
a\leftharpoonup h =(a\triangleleft h)1_A $, for all $a\in A,b\in B$ and $h\in H$.
\end{lemma}
\prf Similarly to \cite{A10}. $\hfill \Box$
\begin{lemma}
Let $H$ be a  Hopf algebra and $A$ an algebra.
Then $(A,\rightharpoonup,\leftharpoonup)$ is  a partial $H$-bimodule
algebra.
\end{lemma}
\prf  Similarly to \cite{A10}. $\hfill \Box$
\medskip

Recall from \cite{A10} that if $B$ is an $H$-module algebra and $A$ is a right ideal of $B$ with
unity $1_A$, the induced partial action on $A$ is called admissible if $B=H\rhd A$.

Similarly, if $B$ is an $H$-module algebra and  $A$ is a left ideal of $B$ with
unity $1_A$, the induced partial action on $A$ is called admissible if $B=A\lhd H$.
\medskip

Now, let $B=H\rhd A$ and $B=A\lhd H$ at the same time, note  it by $B=(A,\triangleright, \triangleleft ).$
\begin{dfn}
Let $H$ be   a Hopf algebra, $B$ an $H$-bimodule  algebra and  $A$ an ideal of $B$ with
unit $1_A.$ The induced partial actions on
$A$ is called admissible if $B=(A,\triangleright, \triangleleft ).$
\end{dfn}
\begin{dfn}
Let $A$ be a partial $H$-bimodule algebra. An enveloping action for $A$ is a pair $(B,\theta),$  where

(a) $B$ is an $H$-bimodule algebra;

(b) The map $\theta:A\rightarrow B$ is a monomorphism of algebras;

(c) The sub-algebra $\theta(A)$ is an ideal in $B$;

(d) The partial action on $A$ is equivalent to the induced partial action on $\theta(A)$;

(e) The induced partial action on $\theta(A)$ is admissible.
\end{dfn}

From now on, we always assume that $( a\leftharpoonup S(h_{(1)}))\o  h_{(2)}=( a\leftharpoonup S(h_{(2)}))\o  h_{(1)}$ for any $a\in A$ and $h\in H$.
This condition  can be easily verified for the case where $H^{\ast}$ is cocommutative (therefore, $H$ is commutative). It is reasonable to assume
this condition beyond the case of $H$ being a cocommutative Hopf algebra.

A concrete counterexample is presented as follows. Recall from \cite{C08} that if $e\in H$ is an central idempotent
such that $e\o e=(e\o 1)\D(e)=\D(e)(e\o 1)$ and $\varepsilon(e)=1$.
A partial $H$-coaction on $A = k$ is given by $\rho^{l}(x)=e\o x\in H\o_k k$ and $\rho^{r}(x)=x\o e\in k\o_k H$.
It is easy to check that  $A$ is a partial $H$-bicomodule algebra.
Therefore, $A$ is a  partial $H^{\ast}$-bimodule algebra. So we can obtain the partial twisted smash product $k \circledast H^{\ast}$, here we only consider the element $f$ of $H^{\ast}$ and check the condition. Then we get
$( a\leftharpoonup S^{\ast}(h_{(2)}))\o  h_{(1)}
= ( x\leftharpoonup S^{\ast}(f_{(2)}))\o  f_{(1)}
=x\o <f_{(2)}, e>f_{(1)}$ and
\begin{eqnarray*}
&&( a\leftharpoonup S^{\ast}(h_{(1)}))\o  h_{(2)}
=( x\leftharpoonup S^{\ast}(f_{(1)}))\o f_{(2)} \\
&&= <f_{(1)}, e>x\o f_{(2)}
= x\o <f_{(1)}, e>f_{(2)}.
\end{eqnarray*}
We can easily obtain $<f_{(1)}, e>f_{(2)}=<f_{(2)}, e>f_{(1)}$ since $e\in H$ is an central idempotent.

\begin{lemma}
Let $\varphi: A\rightarrow Hom(H, A)$ be the map given by $\varphi(a)(h)= h_{(1)}\rightharpoonup a\leftharpoonup S(h_{(2)})$, then we have

(i) $\varphi$ is a linear injective map and an algebra morphism;

(ii) $\varphi(1_A)\ast (h_{(1)}\triangleright \varphi(a)\triangleleft S(h_{(2)}))=\varphi(h_1\rightharpoonup a\leftharpoonup S(h_{(2)}))$ for any $h\in H$ and $a\in A$;

(iii) $\varphi(b)\ast (h_1\triangleright \varphi(a)\triangleleft S(h_{(2)}))=\varphi(b(h_{(1)}\rightharpoonup a\leftharpoonup S(h_{(2)})))$ for any $h\in H$ and $a\in A$.
\end{lemma}

\prf It is easy to see that $\varphi$ is linear, because the partial action is bilinear. Since $\varphi(a)(1_H)=a$, it follows that it is also injective.
For any $a,b\in A$ and $h\in H$, we have
\begin{eqnarray*}
&&\varphi(ab)(h)=h_{(1)}\rightharpoonup(ab)\leftharpoonup S(h_{(2)})\\
&=&[h_{(1)}\rightharpoonup(a)\leftharpoonup S(h_{(4)})][h_{(2)}\rightharpoonup(b)\leftharpoonup S(h_{(3)})]\\
&=&[h_{(1)}\rightharpoonup(a)\leftharpoonup S(h_{(2)})][h_{(3)}\rightharpoonup(b)\leftharpoonup S(h_{(4)})]\\
&=&\varphi(a)(h_{(1)})\varphi(b)(h_{(2)})=\varphi(a)\ast\varphi(b)(h).
\end{eqnarray*}
Therefore $\varphi$ is multiplicative.

For the third claim, we have the following calculation:
\begin{eqnarray*}
&&\varphi(b(h_{(1)}\rightharpoonup a\leftharpoonup S(h_{(2)})))(k)\\
&&=k_{(1)}\rightharpoonup b(h_{(1)}\rightharpoonup a\leftharpoonup S(h_{(2)}))\leftharpoonup S(k_{(2)})\\
&&=[k_{(1)}\rightharpoonup b\leftharpoonup S(k_{(4)})][k_{(2)}\rightharpoonup(h_{(1)}\rightharpoonup a\leftharpoonup S(h_{(2)}))\leftharpoonup S(k_{(3)})]\\
&&=[k_{(1)}\rightharpoonup b\leftharpoonup S(k_{(6)})][k_{(2)}\rightharpoonup 1_A\leftharpoonup S(k_{(5)})][k_{(3)}h_{(1)}\rightharpoonup a\leftharpoonup S(k_{(4)}h_{(2)})]\\
&&=[k_{(1)}\rightharpoonup b\leftharpoonup S(k_{(4)})][k_{(2)}h_{(1)}\rightharpoonup a\leftharpoonup S(k_{(3)}h_{(2)})]\\
&&=[k_{(1)}\rightharpoonup b\leftharpoonup S(k_{(2)})][k_{(3)}h_{(1)}\rightharpoonup a\leftharpoonup S(k_{(4)}h_{(2)})]\\
&&=\varphi(b)(k_1)\varphi(a)(k_{(2)}h)\\
&&=\varphi(b)(k_1)(h\triangleright\varphi(a))(k_{(2)})=\varphi(b)\ast(h\triangleright\varphi(a))(k).
\end{eqnarray*}
Therefore, $\varphi(h\rightharpoonup a)=\varphi(1_A)(h\triangleright \varphi(a))$. One may obtain the second item by setting $b=1_A$.
$\hfill \Box$

\begin{prop}
Let $\varphi: A\rightarrow Hom(H, A)$ be the map defined in Lemma 3.12 and $B=(\varphi(A), \triangleright, \triangleleft)$ the $H$-submodule of $\varphi(A)$. Then

(i) $B$ is an $H$-module subalgebra of $Hom(H, A)$;

(ii) $\varphi(A)$  is a right ideal in $B$ with unity $\varphi(1_A)$.
\end{prop}

\prf Similar to \cite{A10}.
$\hfill \Box$
\medskip

By Lemma 3.12 and Proposition 3.13, we obtain the main result of this section.
\begin{thm}
Let $A$ be a partial $H$-bimodule algebra and $\varphi: A\rightarrow Hom(H, A)$ the map given by $\varphi(a)(h)= h_{(1)}\rightharpoonup a\leftharpoonup S(h_{(2)})$.  Assume that $B=(\varphi(A), \triangleright, \triangleleft)$, then $(B, \varphi)$ is an enveloping
action of $A$.
\end{thm}

\begin{prop}
Let $A$ be a partial $H$-bimodule algebra and $\varphi: A\rightarrow Hom(H, A)$ the map given by $\varphi(a)(h)= h_{(1)}\rightharpoonup a\leftharpoonup S(h_{(2)})$. Assume that $B=(\varphi(A), \triangleright, \triangleleft)$, then $\varphi(A)\unrhd B$ if and only if
\begin{eqnarray*}
k_{(1)}\rightharpoonup (h\rightharpoonup a)\leftharpoonup S(k_{(2)})=[k_{(1)}h_{(1)}\rightharpoonup a\leftharpoonup S(k_{(2)}h_{(2)})][k_{(3)}\rightharpoonup 1_A\leftharpoonup S(k_{(4)})].
\end{eqnarray*}
\end{prop}

\prf Suppose that $\varphi(A)$ is an ideal of B. We know that
\begin{eqnarray*}
\varphi(h\rightharpoonup a)=\varphi(1_A)\ast(h\triangleright \varphi(a))=(h\triangleright \varphi(a))\ast\varphi(1_A).
\end{eqnarray*}
Then, these two functions coincide for all $k\in H$,
\begin{eqnarray*}
\varphi(h\rightharpoonup a)(k)=(h\triangleright \varphi(a))\ast\varphi(1_A)(k).
\end{eqnarray*}
The left-hand side of the previous equality leads to
\begin{eqnarray*}
\varphi(h\rightharpoonup a)(k)=k_{(1)}\rightharpoonup (h\rightharpoonup a)\leftharpoonup S(k_{(2)}).
\end{eqnarray*}
While the right-hand side means
\begin{eqnarray*}
(h\triangleright \varphi(a))\ast\varphi(1_A)(k)&=& (h\triangleright \varphi(a))(k_{(1)})\varphi(1_A)(k_{(2)})\\
&=& \varphi(a)(k_{(1)}h)\varphi(1_A)(k_{(2)})\\
&=& [k_{(1)}h_{(1)}\rightharpoonup a\leftharpoonup S(k_{(2)}h_{(2)})][k_{(3)}\rightharpoonup 1_A\leftharpoonup S(k_{(4)})].
\end{eqnarray*}
Conversely, suppose that the equality
\begin{eqnarray*}
k_{(1)}\rightharpoonup (h\rightharpoonup a)\leftharpoonup S(k_{(2)})=[k_{(1)}h_{(1)}\rightharpoonup a\leftharpoonup S(k_{(2)}h_{(2)})][k_{(3)}\rightharpoonup 1_A\leftharpoonup S(k_{(4)})]
\end{eqnarray*}
holds for all $a\in A$ and $h, k\in H$. Then $\varphi(1_A)$ is a central idempotent in $B$. Therefore $\varphi(A)=\varphi(1_A)B$
is an ideal in $B$.
$\hfill \Box$

\section{A Morita context}
\def\theequation{\arabic{section}.\arabic{equation}}
\setcounter{equation} {0}

In this section, we will construct a Morita context between the partial twisted smash product $\underline{A\circledast H}$ and the twisted smash product $B\circledast H$, where $B$ is an enveloping action
for the partial twisted smash product.
\begin{lemma}
Let $A$ be a partial $H$-bimodule algebra and $(B,\theta)$ an enveloping action, then there is an algebra monomorphism from the
 partial twisted smash product $\underline{A\circledast H}$ into the twisted smash product $B\circledast H.$
 \end{lemma}
\prf Define $\Phi:A\otimes H\rightarrow B\circledast  H$ by $a\otimes h \mapsto \theta(a)\circledast h$ for $h,g\in H$ and $a,b\in A$.
We first check that $\Phi$ is a morphism of algebras as follows:
\begin{eqnarray*}
\Phi((a\otimes h)(b\otimes g))
&=&\Phi( a(h_{(1)}\rightharpoonup b \leftharpoonup S(h_{(3)}))\otimes
h_{(2)}g).  \\
&=& \theta(a(h_{(1)})\rightharpoonup b \leftharpoonup S(h_{(3)}))\circledast
h_{(2)}g\\
&=&  \theta(a)(h_{(1)}\triangleright \theta(b) \triangleleft S(h_{(3)}))\circledast
h_{(2)}g\\
&=& (\theta(a)\circledast h)(\theta(b)\circledast g)\\
&=& \Phi(a\otimes h)\Phi(b\otimes g).
\end{eqnarray*}

Next, we will verify that $\Phi$ is injective. For this purpose, take $x=\sum_{i=1}^{n}a_{i}\otimes h_{i} \in ker\Phi$ and choose $\{a_i\}^{n}_{i=1}$ to be linearly independent. Since $\theta$ is injective, we conclude that $\theta(a_{i})$ are linearly independent.
For each $f\in H^{\ast}$, $\sum_{i=1}^{n}\theta(a_{i})f( h_{i})=0,$  it follows that $f(h_{i})=0,$  so $ h_{i}=0$.
 Therefore we have $x=0$ and $\Phi$ is injective, as desired.

Since the partial twisted smash product $\underline{A\circledast H}$ is a subalgebra of $A\otimes H$, it is
 injectively mapped into $B\circledast H$ by $\Phi.$  A typical element of the image of the partial twisted smash product is
\begin{eqnarray*}
\Phi((a\otimes h)(1_A\otimes 1_{H}))
&=& \Phi(a\otimes h)\Phi(1_A\otimes 1_{H})\\
&=& (\theta(a)\circledast h)(\theta(1_A)\circledast 1_{H})\\
&=& \theta(a(h_{(1)}\triangleright \theta(1_A) \triangleleft S(h_{(3)})))\circledast
h_{(2)}g.
\end{eqnarray*}
And this completes the proof. $\hfill \Box$
\medskip

Take $M=\Phi(A\otimes H)=\{\sum_{i=1}^{n}\theta(a_{i})\circledast  h_{i}; a_{i}\in A\}$  and take $N$ as the subspace of
$B\circledast H$  generated by the elements $(h_{(1)}\triangleright \theta(a) \triangleleft S(h_{(3)}))\otimes
h_{(2)}$ with $h\in H$ and $a\in A.$
\begin{prop}
Let $H$ be a Hopf algebra with an invertible antipode $S$ and $A$ a partial $H$-bimodule algebra.  Suppose that $\theta(A)$ is an ideal of $B$,  then $M$ is a right $B\circledast H$
module and $N$ is a left $B\circledast H$ module.
\end{prop}
\prf In order to prove $M$ is a right $B\circledast H$ module, let $\theta(a)\circledast h\in M$ and $b\circledast k\in B\circledast H$,  then
              $$(\theta(a)\circledast h)(b\circledast k)= \theta(a)(h_{(1)})\triangleright b \triangleleft S(h_{(3)}))\circledast h_{(2)}k. $$
Which lies in $\Phi(A\otimes H)$ because $\theta(A)$ is an ideal in $B$.

Now we show that $N$ is a left $B\circledast H$ module. Let $(h_{(1)}\triangleright \theta(a) \triangleleft S(h_{(3)}))\circledast
h_{(2)}$, where $h\in H$ is  a generator of $N$, then we have
\begin{eqnarray*}
&&(b\circledast k)(h_{(1)}\triangleright \theta(a) \triangleleft S(h_{(3)}))\circledast
h_{(2)})\\
&=& b(k_{(1)}h_{(1)})\triangleright \theta(a) \triangleleft S(k_{(3)}h_{(3)}))\circledast
k_{(2)}h_{(2)}\\
&=& [(\varepsilon(k_{(1)}h_{(1)}\triangleright b)(k_{(2)}h_{(2)})\triangleright \theta(a) ]\triangleleft S(k_{(4)}h_{(4)}))\circledast
k_{(3)}h_{(3)} \\
&=& [((k_{(2)}h_{(2)}S(k_{(1)}h_{(1)})\triangleright b (k_{(3)}h_{(3)})\triangleright\theta(a)]
\triangleleft S(k_{(5)}h_{(5)}))\circledast
k_{(4)}h_{(4)}\\
&=& [(k_{(2)}h_{(2)}\triangleright((S(k_{(1)} h_{(1)})\triangleright b) (k_{(2)}h_{(2)})\triangleright \theta(a)]\triangleleft S(k_{(4)}h_{(4)}))\circledast
k_{(3)}h_{(3)}\\
&=& [(k_{(2)}h_{(2)})\triangleright(S(k_{(1)} h_{(1)})\triangleright b ) \theta(a)]\triangleleft S(k_{(4)}h_{(4)}))\circledast
k_{(3)}h_{(3)}.
\end{eqnarray*}
Because $\theta(A)$ is an ideal of $B$, it follows that $N$ is a left $B\circledast H$ module.
$\hfill \Box$
\medskip

By Proposition 4.2, we can define a left $\underline{A\circledast H}$ module structure on $M$ and a right $\underline{A\circledast H}$ module
structure on N induced by the monomorphism $\Phi$ as follows:
\begin{eqnarray*}
&&(ah_{(1)}\rightharpoonup 1_A\leftharpoonup S(h_{(3)})\otimes h_{(2)})\blacktriangleright (\theta(b))\circledast k)\\
&=&
(\theta(a)h_{(1)}\triangleright \theta(1_A)\triangleleft S(h_{(3)})\circledast h_{(2)})(\theta(b))\circledast k),\\
&&((k_{(1)})\triangleright\theta(b) \triangleleft S(k_{(3)}))\circledast
k_{(2)})\blacktriangleleft( a h_{(1)}\rightharpoonup 1_A\leftharpoonup S(h_{(3)})\otimes h_{(2)})\\
&=&
((k_{(1)})\triangleright \theta(b) \triangleleft S(k_{(3)}))\circledast
k_{(2)})( \theta(a)h_{(1)}\triangleright\theta(1_A)\triangleleft S(h_{(3)})\circledast h_{(2)}).
\end{eqnarray*}
\begin{prop}
Under the same hypotheses of Proposition 4.2, $M$ is indeed a
left $\underline{A\circledast H}$  module with the map $\blacktriangleright$ and N is a right $\underline{A\circledast H}$  module with the map $\blacktriangleleft$.
\end{prop}
\prf We first claim that $\underline{A\circledast H}\blacktriangleright M \subseteq M.$ In fact,
\begin{eqnarray*}
&&(a h_{(1)}\rightharpoonup 1_A\leftharpoonup S(h_{(3)})\otimes h_{(2)})\blacktriangleright (\theta(b))\circledast k)\\
&=&
(\theta(a)(h_{(1)}\triangleright \theta(1_A)\triangleleft S(h_{(3)})\circledast h_{(2)})(\theta(b))\circledast k)\\
&=& \theta(a)(h_{(1)}\triangleright \theta(1_A)\triangleleft S(h_{(5)}))(h_{(2)})\rightharpoonup \theta(b) \leftharpoonup S(h_{(4)}))\circledast h_{(3)}k_{(2)}\\
&=& \theta(a)((h_{(1)})\triangleright \theta(b) \triangleleft S(h_{(3)}))\circledast
h_{(2)}k_{(2)}.
\end{eqnarray*}
Which lies inside $M$ because $\theta(A)$ is an ideal of $B$.
$\hfill \Box$

Next, we verify that $N \blacktriangleleft\underline{A\circledast H}\subseteq N$,  which is similarly to  $N$ is a left $B\circledast H$ module.
which holds because  $\theta(1_A)$ is  a central idempotent.

The last ingredient for a Morita context is to define two bimodule
morphisms
$$\sigma: N\otimes_{\underline{A\circledast H}} M \rightarrow B\circledast H ~~\mbox{and }$$
$$\tau: M\otimes_{B\circledast H} N\rightarrow \underline{A\circledast H}\cong \Phi(\underline{A\circledast H}).$$
As $M, N$ and $\underline{A\circledast H}$ are viewed as subalgebras of $B\circledast H,$  these two maps can be taken
as the usual multiplication on $B\circledast H.$
The associativity of the product assures us
that these maps are bimodule morphisms and satisfy the associativity conditions.
\begin{prop}
 The partial twisted smash
 product $\underline{A\circledast H}$ is Morita equivalent to the  twisted smash product $B\circledast H$.
\end{prop}

\section{Partial representation}
\def\theequation{\arabic{section}.\arabic{equation}}
\setcounter{equation} {0}

In \cite{A10}, Alves and  Batista introduced  the notion of a partial representation of
a Hopf algebra. In this section, we will show a partial representation
of the partial twisted smash product $\underline{A\circledast H}$.

\begin {prop} Let $H$ be a Hopf algebra with an invertible antipode and  $A$  a partial $H$-bimodule algebra. Then the  map
\begin{eqnarray*}
\pi: H\rightarrow End(A),~h\mapsto \pi(h)
\end{eqnarray*}
given by $\pi(h)(a)=h_{(1)}\rightharpoonup a\leftharpoonup S(h_{(2)})$ satisfies:

(1) $\pi(1_H)=1_A;$

(2) $\pi(S^{-1}(h_{(2)}))\pi(h_{(1)})\pi(k)=\pi(S^{-1}(h_{(2)}))\pi(h_{(1)}k).$
\end{prop}
\prf The first identity is quite obvious. In order to prove the second equality, we take $a\in A$ and do the following calculation:
\begin{eqnarray*}
&&\pi(S^{-1}(h_{(2)}))\pi(h_{(1)})\pi(k)(a)\\
&=&S^{-1}(h_{(4)})\rightharpoonup [h_{(1)}\rightharpoonup (k_{(1)}\rightharpoonup a\leftharpoonup S(k_{(2)}))\leftharpoonup S(h_{(2)})]\leftharpoonup h_{(3)}\\
&=& (S^{-1}(h_{(6)})\rightharpoonup 1_A \leftharpoonup h_{(3)})[S^{-1}(h_{(5)})h_{(1)}\rightharpoonup (k_{(1)}\rightharpoonup a\leftharpoonup S(k_{(2)}))\leftharpoonup S(h_{(2)})h_{(4)}]\\
&=& (S^{-1}(h_{(6)})\rightharpoonup 1_A \leftharpoonup h_{(5)})[S^{-1}(h_{(2)})h_{(1)}\rightharpoonup (k_{(1)}\rightharpoonup a\leftharpoonup S(k_{(2)}))\leftharpoonup S(h_{(3)})h_{(4)}]\\
&=& (S^{-1}(h_{(2)})\rightharpoonup 1_A \leftharpoonup h_{(1)})(k_{(1)}\rightharpoonup a\leftharpoonup S(k_{(2)})).
\end{eqnarray*}
On the other hand, we have
\begin{eqnarray*}
&&\pi(S^{-1}(h_{(2)}))\pi(h_{(1)}k)(a)\\
&&=S^{-1}(h_{(4)})\rightharpoonup [h_{(1)}k_{(1)}\rightharpoonup a\leftharpoonup S(h_{(2)}k_{(2)}))]\leftharpoonup h_{(3)}\\
&&=(S^{-1}(h_{(6)})\rightharpoonup 1_A \leftharpoonup h_{(3)})[S^{-1}(h_{(5)})h_{(1)}k_{(1)}\rightharpoonup a\leftharpoonup S(h_{(2)}k_{(2)})h_{(4)}]\\
&&=(S^{-1}(h_{(6)})\rightharpoonup 1_A \leftharpoonup h_{(5)})[S^{-1}(h_{(2)})h_{(1)}k_{(1)}\rightharpoonup a\leftharpoonup S(h_{(3)}k_{(2)})h_{(4)}]\\
&&=(S^{-1}(h_{(2)})\rightharpoonup 1_A \leftharpoonup h_{(1)})(k_{(1)}\rightharpoonup a\leftharpoonup S(k_{(2)})).
\end{eqnarray*}
And this completes the proof. $\hfill \Box$
\medskip

With this result we can propose the following definition  similar to \cite{A10}.
\begin{dfn} Let $H$ be a Hopf algebra with an invertible antipode. A partial
representation of $H$ on a unital algebra $B$ is a linear map
\begin{eqnarray*}
\pi: H\rightarrow B,~h\mapsto \pi(h)
\end{eqnarray*}
such that

(1) $\pi(1_H)=1_B;$

(2) $\pi(S^{-1}(h_{(2)}))\pi(h_{(1)})\pi(k)=\pi(S^{-1}(h_{(2)}))\pi(h_{(1)}k).$
\end{dfn}
\begin {thm} Let $H$ be a Hopf algebra with an invertible antipode and  $A$  a partial $H$-bimodule algebra. Then the linear map
\begin{eqnarray*}
\pi: H\rightarrow \underline{A\circledast H},~
h\mapsto (h_{(1)}\rightharpoonup 1_A\leftharpoonup S(h_3))\o h_{(2)}
\end{eqnarray*}
is a partial representation of $H$.
\end{thm}
\prf The first identity is quite obvious. In order to prove the second equality, we need the following calculations:
\begin{eqnarray*}
&&\pi(S^{-1}(h_{(2)}))\pi(h_{(1)})\pi(k)\\
&&=((S^{-1}(h_{(6)})\rightharpoonup 1_A\leftharpoonup h_{(4)})\o S^{-1}(h_{(5)}))(h_{(1)}\rightharpoonup 1_A\leftharpoonup S(h_{(3)})\o h_{(2)})\\
&&\hspace{9cm}(k_{(1)}\rightharpoonup 1_A\leftharpoonup S(k_{(3)})\o k_{(2)})\\
&&=(S^{-1}(h_{(8)})\rightharpoonup 1_A\leftharpoonup h_{(4)})[S^{-1}(h_{(7)})\rightharpoonup(h_{(1)}\rightharpoonup 1_A\leftharpoonup S(h_{(3)})\leftharpoonup h_{(5)}]\\
&&\hspace{6cm}\o S^{-1}(h_{(6)})h_{(2)})(k_{(1)}\rightharpoonup 1_A\leftharpoonup S(k_{(3)})\o k_{(2)})\\
&&=(S^{-1}(h_{(10)})\rightharpoonup 1_A\leftharpoonup h_{(4)})(S^{-1}(h_{(9)})\rightharpoonup 1_A\leftharpoonup h_{(5)})
\end{eqnarray*}
\begin{eqnarray*}
&&[(S^{-1}(h_{(8)})h_{(1)}\rightharpoonup 1_A\leftharpoonup S(h_{(3)})h_{(6)}]\o S^{-1}(h_{(7)})h_{(2)})(k_{(1)}\rightharpoonup 1_A\leftharpoonup S(k_{(3)})\o k_{(2)})\\
&&=(S^{-1}(h_{(8)})\rightharpoonup 1_A\leftharpoonup h_{(4)})[(S^{-1}(h_{(7)})h_{(1)}\rightharpoonup 1_A\leftharpoonup S(h_{(3)})h_{(5)}]\\
&&\hspace{6cm}\o S^{-1}(h_{(6)})h_{(2)})(k_{(1)}\rightharpoonup 1_A\leftharpoonup S(k_{(3)})\o k_{(2)})\\
&&=(S^{-1}(h_{(7)})\rightharpoonup 1_A\leftharpoonup h_{(2)})[(S^{-1}(h_{(6)})h_{(1)}\rightharpoonup 1_A\leftharpoonup S(h_{(3)})h_{(4)}]\\
&&\hspace{6cm}\o S^{-1}(h_{(6)})h_{(5)})(k_{(1)}\rightharpoonup 1_A\leftharpoonup S(k_{(3)})\o k_{(2)})\\
&&=(S^{-1}(h_{(2)})\rightharpoonup 1_A\leftharpoonup h_{(1)})(k_{(1)}\rightharpoonup 1_A\leftharpoonup S(k_{(3)})\o k_{(2)}).
\end{eqnarray*}
On the other hand, we have
\begin{eqnarray*}
&&\pi(S^{-1}(h_{(2)}))\pi(h_{(1)}k)\\
&&=[S^{-1}(h_{(6)})\rightharpoonup 1_A\leftharpoonup h_{(4)}\o S^{-1}(h_{(5)})][h_{(1)}k_{(1)}\rightharpoonup 1_A\leftharpoonup S(h_{(3)}k_{(3)})\o h_{(2)}k_{(2)}]\\
&&=(S^{-1}(h_{(8)})\rightharpoonup 1_A\leftharpoonup h_{(4)}) S^{-1}(h_{(7)})\rightharpoonup[h_{(1)}k_{(1)}\rightharpoonup 1_A\leftharpoonup S(h_{(3)}k_{(3)})]\leftharpoonup h_{(5)}\\
 &&\hspace{10cm}   \o S^{-1}(h_{(6)})h_{(2)}k_{(2)}]\\
&&=(S^{-1}(h_{(10)})\rightharpoonup 1_A\leftharpoonup h_{(4)})(S^{-1}(h_{(9)})\rightharpoonup 1_A\leftharpoonup h_{(5)}) \\
 &&\hspace{3cm}  [S^{-1}(h_{(8)})h_{(1)}k_{(1)}\rightharpoonup 1_A\leftharpoonup S(h_{(3)}k_{(3)})h_{(6)}]\o S^{-1}(h_{(7)})h_{(2)}k_{(2)}]\\
&&=(S^{-1}(h_{(8)})\rightharpoonup 1_A\leftharpoonup h_{(4)})[S^{-1}(h_{(7)})h_{(1)}k_{(1)}\rightharpoonup 1_A\leftharpoonup S(h_{(3)}k_{(3)})h_{(5)}]\\
 &&\hspace{10cm}  \o S^{-1}(h_{(6)})h_{(2)}k_{(2)}]\\
 &&=(S^{-1}(h_{(8)})\rightharpoonup 1_A\leftharpoonup h_{(7)})[S^{-1}(h_{(2)})h_{(1)}k_{(1)}\rightharpoonup 1_A\leftharpoonup S(h_{(3)}k_{(3)})h_{(4)}]\\
 &&\hspace{10cm}  \o S^{-1}(h_{(6)})h_{(5)}k_{(2)}]\\
 &&=(S^{-1}(h_{(2)})\rightharpoonup 1_A\leftharpoonup h_{(1)})(k_{(1)}\rightharpoonup 1_A\leftharpoonup S(k_{(3)})\o k_{(2)}).
\end{eqnarray*}
So $\pi(S^{-1}(h_{(2)}))\pi(h_{(1)})\pi(k)=\pi(S^{-1}(h_{(2)}))\pi(h_{(1)}k)$, as required. $\hfill \Box$

\section{Duality for partial twisted smash products}
\def\theequation{\arabic{section}.\arabic{equation}}
\setcounter{equation} {0}

In this section, we explore some consequences
of globalization theorems in order to present versions of the duality theorems
of  Blattner-Montgomery for partial twisted smash products, generalizing the results of \cite{L08}.

Let $H$ be a Hopf algebra which is finitely generated and projective as $k$-module with
dual basis $\{(b_i,p_i)\in H\o H^\ast|1\leq i\leq n\}$. Assume that $H^\ast$ acts on $H$ from the left by $f\rightarrow h=\sum h_{(1)}f(h_{(2)})$ and the right by $ h \leftarrow f=\sum h_{(2)}f(h_{(1)})$,  such that the smash product $H\# H^\ast$ can be considered as an algebra whose multiplication is given by
\begin{eqnarray*}
(h\# f)(k\# g)=\sum h(f_{(1)}\rightarrow k)\# f_{(2)}\ast g,
\end{eqnarray*}
for any $h,k\in H$ and $f, g\in H^\ast$.

\begin{lemma}\cite{S69}
Let $H$ be a finite dimensional Hopf algebra. Then the linear maps

(1) $\lambda: H\# H^\ast \rightarrow End(H), ~\lambda(h\# f)(k)=h(f\rightarrow k),$

(2) $\varphi: H^\ast\# H \rightarrow End(H), ~\varphi(f\# h)(k)=(k\leftarrow f)h,$

\noindent are isomorphisms of algebras, where $h,k\in H$ and $f, g\in H^\ast$.
\end{lemma}
 The partial twisted smash product $ \underline{A\circledast H}$ in Proposition 3.3. becomes naturally a right $H$-comodule algebra by
\begin{eqnarray*}
\rho=1\o \D: A\circledast H\o H\rightarrow A\o H\o H, ~~~a\o h\mapsto a\o h_{(1)}\o h_{(2)}.
\end{eqnarray*}
 For $(a\o h)1_A\in \underline{A\circledast H}$, we have
\begin{eqnarray*}
\rho((a\o h)1_A)= a(h_{(1)}\rightharpoonup 1_A\leftharpoonup S(h_{(2)}))\o h_{(3)}\o h_{(4)},
\end{eqnarray*}
which make $\underline{A\circledast H}$ into a right $H$-comodule algebra.
 Moreover, $\underline{A\circledast H}$ becomes a left $H^\ast$ module algebra, where the action is defined by
\begin{eqnarray*}
f\c((a\o h)1_A)=a(h_{(1)}\rightharpoonup 1_A\leftharpoonup S(h_{(3)}))\#(f\rightarrow h_{(2)})=(a\#(f\rightarrow h))1_A.
\end{eqnarray*}
for all $f\in H^\ast, h\in H, a\in A$.

Similar to \cite{L08},  we can define a homomorphism $\phi: A\rightarrow A\o End(H)$ by
\begin{eqnarray*}
\phi(a)=\sum_{i=1}^{n} (b_{i(1)}\rightharpoonup a\leftharpoonup S(b_{i(2)}))\o \varphi(S^{-1}(p_i)\o 1_H).
\end{eqnarray*}
Then $\phi$ is an algebra homomorphism.
\begin{lemma}
Let $\psi: H\#H^{\ast}\rightarrow A\o End(H)$ be the map defined by $h\#f \mapsto 1\o \lambda(h\# f)$ for all $h\in H$ and $f\in H^\ast$. Then we have
\begin{eqnarray*}
\phi(1_A)\psi(h\#f)\phi(a)=\phi(h_{(1)}\rightharpoonup a\leftharpoonup S(h_{(3)}))\psi(h_{(2)}\#f).
\end{eqnarray*}
\end{lemma}
\prf For any $h\in H, f\in H^\ast$ and $a\in A$, we have
\begin{eqnarray*}
&&\phi(h_{(1)}\rightharpoonup a\leftharpoonup S(h_{(3)}))\psi(h_{(2)}\#f)\\
&=&\sum_{i} p_i(h_{(1)})\phi(b_{i(1)}\rightharpoonup a\leftharpoonup S(b_{i(3)}))\psi(h_{(2)}\#f)\\
&=& \sum_{i,j}b_{j(1)}\rightharpoonup(b_{i(1)}\rightharpoonup a\leftharpoonup S(b_{i(3)}))\leftharpoonup S(b_{j(3)})\o \varphi(S^{-1}(p_j)\#1_H)\lambda(h\leftarrow p_i\# f)\\
&=& \sum_{k,r}(b_{k(1)}\rightharpoonup1_A \leftharpoonup S(b_{k(2)})) (b_{r(1)}\rightharpoonup a\leftharpoonup S(b_{r(2)}))\o \varphi(S^{-1}(p_k)\#1_H)\\
&&\hspace{7cm}\varphi(S^{-1}(p_{k(1)})\#1_H)\lambda(h\leftarrow p_{k(2)}\# f)\\
&=&\phi(1_A) \sum_{r}(b_{r(1)}\rightharpoonup a\leftharpoonup S(b_{r(2)}))\o \varphi(S^{-1}(p_{r})_{(2)}\#1_H)\lambda(h\leftarrow S(S^{-1}(p_{r})_{(1)})\# f)\\
&=&\phi(1_A) \sum_{r}(b_{r(1)}\rightharpoonup a\leftharpoonup S(b_{r(2)}))\o \lambda(h\# f)\varphi(S^{-1}(p_{r})\#1_H)\\
&=& \phi(1_A)\psi(h\#f)\phi(a).
\end{eqnarray*}
The proof is completed.  $\hfill \Box$

\begin{thm}
Let $H$ be a finitely generated projective  Hopf algebra and $A$  a partial $H$-bimodule algebra.
Then the map
 $$\Phi: A\o H\# H^\ast\rightarrow A\o End(H),~
                      a\o h\# f\mapsto \phi(a)\psi(h\#f)$$
is an algebra homomorphism.
 The image of the restriction to $\underline{A\circledast H}\# H^{\ast}$ lies inside $e(A\o End(H))e$, where $e$ is the idempotent defined by
\begin{eqnarray*}
e=\sum_{i=1}^{n}(b_{i(1)}\rightharpoonup 1_A\leftharpoonup S(b_{i(2)}))\o \varphi(S^{-1}(p_i)\o 1_A).
\end{eqnarray*}
\end{thm}
\prf For any $a,b\in A, h,k\in H$ and $f, g\in H^\ast$, we have
\begin{eqnarray*}
\Phi(a\o h\# f)\Phi(b\o k\# g)&=& \phi(a)\psi(h\# f)\phi(b)\psi(k\# g)\\
&=& \phi(a) \phi(1_A)\psi(h\# f)\phi(b)\psi(k\# g)\\
&=& \phi(a) \phi(h_{(1)}\rightharpoonup b\leftharpoonup S(h_{(3)}))\psi(h_{(2)}\#f)\psi(k\# g)\\
&=& \phi(a(h_{(1)}\rightharpoonup b\leftharpoonup S(h_{(3)})))\psi(h_{(2)}(f_{(1)}\rightarrow k)\# f_{(2)}\ast g)\\
&=& \Phi(\phi(a(h_{(1)}\rightharpoonup b\leftharpoonup S(h_{(3)})))\o h_{(2)}(f_{(1)}\rightarrow k)\# f_{(2)}\ast g)\\
&=& \Phi((a\o h\# f)(b\o k\# g)).
\end{eqnarray*}
Hence $\Phi$ is an algebra homomorphism.
Since the image of the identity $1=1_A\circledast 1_H\# 1_{H^\ast}$ of $\underline{A\circledast H}\# H^{\ast}$ under the map $\Phi$ is $e$,
$e$ is an idempotent.
 Moreover, for any $\g\in \underline{A\circledast H}\# H^{\ast}$, we have
\begin{eqnarray*}
\Phi(\gamma)=\Phi(1_A\gamma1_A)\in e(A\o End(H))e,
\end{eqnarray*}
as desired. And this completes the proof.$\hfill \Box$

\section*{\bf Acknowledgement}

 The authors would like to thank Professor Ja Kyung Koo for the helpful
 comments. The work was  supported  by the NSF of Jiangsu Province (BK2012736) and the Fund of Science and Technology Department of Guizhou Province (No. 2014GZ81365).

\end{document}